\documentclass[11pt, a4paper]{amsart}
\usepackage{amssymb, amsmath, amscd, amsthm, amsxtra}

\newtheorem{theorem}{Theorem}[section]

\numberwithin{equation}{subsection}

\begin{document}

\title{A Remark on Soliton Equation of Mean Curvature Flow}
\author{Li Ma  and Y.Yang}
\address{Department of Mathematics\\Tsinghua
University\\Beijing,100084,China}
\email{lma@math.tsinghua.edu.cn}
\date{}

\thanks{The work of Ma is partially supported by the Key 973 project of China}

\keywords{Self-Similar, Mean curvature flow}
\subjclass{53C44,
53C42} \maketitle

\begin{abstract}
In this short note, we consider self-similar immersions $F:
\mathbb{R}^n \to \mathbb{R}^{n+k}$ of the Graphic Mean Curvature
Flow of higher co-dimension. We show that the following is true:
Let $F(x) = (x,f(x)), x \in \mathbb{R}^{n}$ be a graph solution to
the soliton equation
$$
\overline{H}(x) + F^{\bot}(x) = 0.
$$
Assume $\sup_{\mathbb{R}^{n}}|Df(x)| \le C_{0} < + \infty$. Then
there exists a unique smooth function $f_{\infty}:
\mathbb{R}^{n}\to \mathbb{R}^k$ such that
$$
f_{\infty}(x) = \lim_{\lambda \to \infty}f_{\lambda}(x)
$$
and
$$
f_{\infty}(r x)=r f_{\infty}(x)
$$
for any real number $r\not= 0$, where
$$
f_{\lambda}(x) = \lambda^{-1}f(\lambda x).
$$
\end{abstract}

\section{Introduction}
Let $M^{n+k}$ be a Riemannian manifold of dimension $n+k$. Assume
that $\Sigma^{n}$ be a Riemannian manifold  of dimension $n$
without boundary. Let $F:\Sigma^{n} \to M^{n+k}$ be an isometric
immersion. Denote $\nabla$ (respectively $D$) the covariant
differentiation on $\Sigma$ (on  $M$). Let $T\Sigma$ and $N\Sigma$
be the tangent bundle and normal bundle of $\Sigma$ in $M$
respectively. We define the second fundamental form of the
immersion $\Sigma$ by
$$
II: T\Sigma \times T\Sigma \to N\Sigma, $$
with
$$II(X,Y) = D_{X}Y - \nabla_{X}Y,
$$
for tangential vector fields $X,Y$ on $\Sigma$. We define the mean
curvature vector field (in short, MCV) by
$$
\overline{H} = \mbox{tr}_{\Sigma}II.
$$

In recent years, many people are interested in studying the
evolution of the immersion $F:\Sigma^{n} \to M^{n+k}$ along its
Mean Curvature Flow (in short, just say MCF).
 The MCF is
defined as follows. Given an one-parameter family of sub-manifolds
$\Sigma_t=F_t(\Sigma)$ with immersions $F_t:\Sigma \longrightarrow
M$ . Let $\overline{H}(t)$ be the MCV of $\Sigma_t$. Then our MCF
is the equation/system
\[
\frac{\partial F(x,t)}{\partial t}=\overline{H}(x,t).
\]
This flow has many very nice results if the codimension $k=1$. See
the work of G.Huisken \cite{H2} for a survey in this regard. Since
there is very few result about MCF in higher codimension, we will
study it in the target when $M^{n+k}=R^{n+k}$, which is the
standard Euclidian space.

In this short note, we will consider a family of self-similar
graphic immersions $F(\cdot, t): \mathbb{R}^n \to
\mathbb{R}^{n+k}$ of the Mean Curvature Flow (MCF):
$$
 \frac{\partial }{\partial t}F(x,t) =
\overline{H}(x,t),\;\; \forall x \in \mathbb{R}^n,\;\; \forall t
\in (-\infty,0).
$$
Write
$$
\Sigma_{t} = F(\mathbb{R}^n,t), $$ and
$$F = (F^{A}), \;\;\; 1 \le A \le n+k.
$$
By definition, we call the family $\Sigma_{t}$ {\it self-similar}
if
$$
\Sigma_{t} = \sqrt{-t} \Sigma_{-1}, \;\;\;\forall t < 0.
$$
In this case, we can reduce the MCF into an elliptic system. In
the other word, we have the following parametric elliptic equation
for the family $\Sigma_{t}$:
$$
\overline{H}(x) + F^{\bot}(x) = 0, \;\;\;\forall x \in \Sigma
_{-1} := \Sigma.
$$
We will call this system as the {\it soliton equation} of the MCF.
Note that this equation is usually obtained from the monotonicity
formula of G.Huisken \cite{H1} for blow-up. It is a hard and open
problem to classify solutions of this equation.

Fix $\Sigma=\Sigma_t$. Assume that $F(x)=(x,f(x))$. Let
$$
Q=(Q^{A}_{\alpha}),\;\;{n+1\leq\alpha\leq n+k}\;\;1\leq A\leq n+k
$$ is the orthogonal
projection onto $N_{p}\Sigma$, where $p\in \Sigma$. Then the
second fundamental form of $\Sigma$ can be written as
$$
\mbox{II}_{ij}^{A} = Q_{\alpha}^{A}D_{ij}^{2}f^{\alpha}.
$$
Hence, we have the expression for the mean curvature vector of
$\Sigma$ in $\mathbb{R}^{n+k}$:
$$
\overline{H}^{A} = g^{ij}Q^{A}_{\alpha}D^{2}_{ij}f^{\alpha}.
$$

Our main result in this paper is the following
\begin{theorem}
Let $F(x) = (x,f(x)), x \in \mathbb{R}^{n}$ be a graph solution to
the soliton equation
$$
\overline{H}(x) + F^{\bot}(x) = 0.
$$
Assume $\sup_{\mathbb{R}^{n}}|Df(x)| \le C_{0} < + \infty$. Then
there exists a unique smooth function $f_{\infty}:
\mathbb{R}^{n}\to \mathbb{R}^k$ such that
$$
f_{\infty}(x) = \lim_{\lambda \to \infty}f_{\lambda}(x)
$$
and
$$
f_{\infty}(r x)=r f_{\infty}(x)
$$
for any real number $r\not= 0$, where
$$
f_{\lambda}(x) = \lambda^{-1}f(\lambda x).
$$
\end{theorem}

We remark that the proof of this result given below is very
simple. But it is based on a nice observation. We just use the
divergence theorem with a nice test function. In the next section,
we recall the form of divergence theorem for convenient of the
readers. In the last section we give a proof of our Theorem.

We point out that we may consider
$F_{\infty}(x)=(x,f_{\infty}(x))$ obtained above as a tangential
minimal cone along the research direction done by L.Simon
\cite{Si}.

\section{Preliminary}
Given a vector field $X: \Sigma \to TM$. Let $X^{T}$ and $X^{N}$
denote the projection of $X$ onto $T\Sigma$ and $N\Sigma$ respectively.
We define the divergence of $X$ on $\Sigma$ as
$$
\mbox{div}_{\Sigma}X = \sum g^{ij}\langle D_{i}X, \frac{\partial }{\partial x^{j}}\rangle
$$
where $(g^{ij}) = g_{ij}^{-1}$, and $(g^{ij})$ is the induced metric
tensor written in local coordinates $(x^{i})$ on $\Sigma$.

Note that, for any tangential vector field $Y$ on $\Sigma$,
$$
D_{Y}X = D_{Y}X^{T} + D_{Y}X^{N}.
$$
So
\begin{align*}
\langle D_{Y}X,Y\rangle & = \langle D_{Y}X^{T},Y \rangle + \langle D_{Y}X^{N},Y \rangle \\
& = \langle \nabla_{Y}X^{T},Y \rangle - \langle D_{Y}Y, X^{N} \rangle \\
& = \langle \nabla_{Y}X^{T},Y \rangle - \langle II(Y,Y),X \rangle.
\end{align*}
Hence
$$
\mbox{div}_{\Sigma}X^{T} = \mbox{div}_{\Sigma}X + \langle
X,\overline{H} \rangle,
$$
and by the Stokes formula on $\Sigma$, we have
$$
\int_{\Sigma}\mbox{div}X^{T} = \int_{\partial_{\Sigma}}\langle X,\nu\rangle d\sigma
$$
and
$$
\int_{\Sigma}\mbox{div}_{\Sigma}X d\nu = - \int_{\Sigma}\langle
\overline{H},X\rangle d\nu + \int_{\partial \Sigma}\langle
X,\nu\rangle d\sigma,
$$
where $\nu$ is the exterior normal vector field to $\Sigma$ on
$\partial\Sigma$.

\section{Proof of Main Theorem}

In the following, we take $M^{n+k} = \mathbb{R}^{n+k}$ as the
standard Euclidean space. We assume that the assumption of our
Theorem 1.1 is true in this section.

Define the vector field
$$
X = (1 + |F|)^{-s}F
$$
where $s \in \mathbb{R}$ to be determined.

Note that, $\nabla |F| = \frac{F^{\top}}{|F|}$ and $\mbox{div}_{\Sigma}F
= n$.
So
\begin{align*}
\mbox{div}_{\Sigma}X & = \langle\nabla(1+|F|)^{-s}, F\rangle + (1+|F|)^{-s}
\mbox{div}_{\Sigma}F \\
& = - \frac{s(1+|F|)^{-s-1}}{|F|}|F^{\top}|^{2} + n(1 + |F|)^{-s}.
\end{align*}
Locally, we may assume that $\Sigma$ is a graph of the form
$(x,f(x)) \in B_{R}(0) \times \mathbb{R}^{k}$, where $B_{R}(0)$ is
the ball of radius $R$ centered at $0$. Let $\Sigma_{R} = \Sigma
\cap (B_{R}(0) \times \mathbb{R}^{k})$. By the divergence theorem
we have
$$
\int_{\Sigma_{R}} \mbox{div}_{\Sigma}X  =
\int_{\Sigma_{R}}\langle\overline{H},X\rangle - \int_{\partial
\Sigma_{R}}\langle X,\nu\rangle
$$
By direct computation, we have that
\begin{align*}
\int_{\Sigma_{R}} \mbox{div}_{\Sigma}X & = -s
\int_{\Sigma_{R}}\frac{(1+|F|)^{-s-1}}{|F|}|F^{\top}|^{2}
  + n \int_{\Sigma_{R}}(1 + |F|)^{-s} \\
& = - \int_{\Sigma_{R}}(1+|F|)^{-s}|F^{\bot}|^{2}
  - \int_{\partial_{\Sigma_{R}}}(1+|F|)^{-s}\langle F,\nu\rangle \\
& = - \int_{\Sigma_{R}}(1+|F|)^{-s}|\overline{H}|^{2}
  - \int_{\partial_{\Sigma_{R}}}(1+|F|)^{-s}\langle F,\nu\rangle.
\end{align*}
Hence, we have
$$
\int_{\Sigma_{R}}(1+|F|)^{-s}|\overline{H}|^{2} =
  s \int_{\Sigma_{R}}\frac{(1+|F|)^{-s-1}}{|F|}|F^{\top}|^{2}
  - n \int_{\Sigma_{R}}(1 + |F|)^{-s}
  - \int_{\partial_{\Sigma_{R}}}(1+|F|)^{-s} \langle F,\nu\rangle.
$$
Since $|F^{\top}| \le |F| \le 1+|F|$, we have
$$
\int_{\Sigma_{R}}\frac{(1+|F|)^{-s-1}}{|F|}|F^{\top}|^{2} \le
\int_{\Sigma_{R}}(1 + |F|)^{-s}.
$$
Clearly we have
$$
\left| \int_{\partial_{\Sigma_{R}}}(1+|F|)^{-s} \langle
F,\nu\rangle \right| \le
\int_{\partial_{\Sigma_{R}}}(1+|F|)^{1-s}.
$$
Combining these two inequalities together we get
$$
\int_{\Sigma_{R}}(1+|F|)^{-s}|\overline{H}|^{2} \le
(s-n) \int_{\Sigma_{R}}(1 + |F|)^{-s}
+ \int_{\partial_{\Sigma_{R}}}(1+|F|)^{1-s}.
$$
Choosing $s=n$ yields $(*)$:
$$
 \int_{\Sigma_{R}}(1+|F|)^{-n}|\overline{H}|^{2} \le
\int_{\partial_{\Sigma_{R}}}(1+|F|)^{1-n}.
$$
By our assumption we have that $\exists C > 0$ such that for $F(x)
= (x,f(x))$ on $\Sigma = \mathbb{R}^{n}$, we have $$\det(I +
(df)^{\top}df) \le C$$ on $\Sigma$. Since
$$
g_{ij} = \delta_{ij} + D_{i}f^{\alpha} \cdot D_{j}f^{\alpha},
$$
we know that
$$
I \le (g_{ij}) \le CI.
$$
Hence
$$
(1 + |x|) \le (1 + |F(x)|) \le C(1+|x|).
$$
Therefore we get from $(*)$ the key estimate $(K)$:
$$
\int_{B_{R}(0)}(1+|x|)^{-n}|\overline{H}|^{2}dx \le C
\int_{\partial{B}_{R}(0)}(1+|x|)^{1-n} \le C.
$$

We now go to the proof of our Theorem.
\begin{proof}
Note that the mean curvature flow for the graph of $f$ can be read
as
$$
\frac{\partial f^{\alpha}}{\partial t} =
g^{ij}D_{ij}^{2}f^{\alpha}, \alpha = 1, \cdots, k.
$$
The important fact about this equation is that it is invariant
under the transformation
$$
f(x) \to \frac{1}{\lambda}f(\lambda x), \forall \lambda > 0.
$$
Compute
\begin{align*}
\frac{d}{d \lambda}f_{\lambda}(x) & = - \lambda^{-2}f(\lambda x)
  + \lambda^{-1} Df(\lambda x)\cdot x \\
  & = \lambda^{-2} [Df(\lambda x) \cdot \lambda x - f(\lambda x)] \\
  & = \lambda^{-2} \langle(Df(\lambda x),-1),(\lambda x, f(\lambda x))
    \rangle \\
  & = \lambda^{-2} \langle(Df(\lambda x),-1),F(\lambda x) \rangle \\
  & = \lambda^{-2} \langle(Df(\lambda x),-1),F(\lambda x)^{\bot}
    \rangle.
\end{align*}
Here we have used the fact that
$$
(Df(\lambda x),-1) \bot T_{p} \Sigma.
$$
So
$$
\frac{d}{d \lambda}f_{\lambda}(x) = \lambda^{-2}\langle(-Df(\lambda x),
1), \overline{H} \rangle.
$$
Hence
$$
\left| \frac{d}{d \lambda}f_{\lambda}(x) \right| \le C\lambda^{-2}
|\overline{H}|.
$$
So, for $x \in S^{n-1}$, we have
\begin{align*}
|f_{\lambda}(x) - f_{\mu}(x)| & \le C \int^{\mu}_{\lambda}\frac{
  \overline{H}(\lambda x)}{\sigma^{2}} d\sigma \\
  & \le C(\int_{\lambda}^{\mu} \frac{1}{\sigma^{3}}d\sigma)
    (\int_{\lambda}^{\mu}\frac{|\overline{H}^{2}|(\sigma x)}{\sigma}
    d\sigma) \\
  & \le C|\mu^{-2} - \lambda^{-2}|\int_{\lambda}^{\mu}\frac{|
    \overline{H}(\sigma x)|^{2}}{\sigma} d\sigma.
\end{align*}
Notice that, for $\mu \ge \lambda > 1$,
$$
\int_{S^{n-1}}dx\int_{\lambda}^{\mu}\frac{|\overline{H}(\sigma
x)|^{2}}{\sigma}
  d\sigma \le \int_{0}^{\infty}\int_{S^{n-1}}\frac{|\overline{H}
  (\sigma x)|^{2}}{(1 + \sigma)^{n}} \sigma^{n-1} dx d\sigma \le C.
$$
The last inequality follows from the inequality (K). Therefore, we
have the estimate $(**)$:
$$
\int_{S^{n-1}}|f_{\lambda}(x) - f_{\mu}(x)|^2 dx\leq C|\mu^{-2} -
\lambda^{-2}|.
$$

This implies that $(f_{\lambda})$ is a Cauchy sequence in $L^{2}
(S^{n-1})$. Let $f_{\infty}$ be its unique limit. Since
$\sup_{\mathbb{R}^{n}} |Df_{\lambda}| = \sup_{\mathbb{R}^{n}}|Df|
\le C_{0}$, the Arzela-Ascoli theorem tells us that
$(f_{\lambda})$ is compact in $C^{\alpha}(S^{n-1}), \forall \alpha
\in (0,1)$. Therefore
$$
f_{\infty}(x) = \lim f_{\lambda}(x)\;\;\; \mbox{uniformly on }
S^{n-1},
$$
and
$$
f_{\infty}(rx) = r f_{\infty}(x), \;\;\;\forall 0\not=r \in
\mathbb{R}.
$$
This finishes the proof of Theorem 1.1
\end{proof}

In the following, we pose a question about the stability of
self-similar solutions of (MCF). Let $f_{0}:\mathbb{R}^{n} \to
\mathbb{R}^{k}$ be a smooth function with uniformly bounded
(Lipschitz) gradient. Assume
$$
\lim_{\lambda \to \infty} f_{0\lambda} = f^{\infty}_{0},\;\;\;
\mbox{uniformly on } S^{n-1}.
$$
Assume $f:\mathbb{R}^{n} \times [0, \infty) \to \mathbb{R}^{k}$ such
that $F(x,t) = (x,f(x,t))$ is a solution of (MCF) with the initial data
$F(x,0) = (x,f_{0}(x)).$ We ask if there is a smooth mapping $\hat
{f}: \mathbb{R}^{n} \to \mathbb{R}^{k}$ such that $\hat{f}(\cdot,s)
\to \hat{f}(\cdot)$ uniformly on compact subsets of $\mathbb{R}^{n}$
as $s \to \infty$. Here $\hat{f}$ is defined by
$$
\hat{f}(x,s) = t^{-\frac{1}{2}}f(\sqrt{t} x, t), s = \frac{1}{2}
\log t, 0 \le s < \infty\;\; \mbox{with }\;\; t \ge 1.
$$
A related stability result is done by one of us in \cite{M}.

\end{document}